# THE BERNSTEIN–VON MISES THEOREM FOR THE PROPORTIONAL HAZARD MODEL[1]

By Yongdai Kim

*Seoul National University*

We study large sample properties of Bayesian analysis of the proportional hazard model with neutral to the right process priors on the baseline hazard function. We show that the posterior distribution of the baseline cumulative hazard function and regression coefficients centered at the maximum likelihood estimator is jointly asymptotically equivalent to the sampling distribution of the maximum likelihood estimator.

**1. Introduction.** Since Cox [3] proposed the proportional hazard model for survival time data in the presence of covariates, the proportional hazard model has enjoyed a wide variety of applications in biomedical data analysis and reliability. Although it does not require any parametric assumption on the baseline cumulative hazard function (c.h.f.), its computation is almost parametric. By casting the theoretical framework as a counting process problem, the study of its asymptotic properties becomes a historical success story in theoretical statistics. These are some of many reasons for its popularity in applications as well as the theory of statistics.

The Bayesian analysis of the proportional hazard model has also been studied by many authors. Kalbfleisch [11] studied its Bayesian analysis with gamma process priors on the baseline c.h.f. For the Bayesian analysis of the proportional hazard model with beta process priors, a Markov chain Monte Carlo computation is proposed by Laud, Damien and Smith [17] and Lee and Kim [18], and the marginal posterior distribution of the regression coefficients is obtained by Hjort [10]. Kim and Lee [14] obtained the posterior

Received April 2004; revised March 2005.

[1]Supported in part by the SRC/ERC program of MOST/KOSEF (R11-2000-073-00000) and in part by Grant R01-2004-000-10284-0 from the Basic Research Program of the Korea Science and Engineering Foundation.

distribution for the proportional hazard model with neutral to the right process priors [5] when the survival times are under left truncation and right censoring. Most research effort from the Bayesian side has been devoted to identifying the posterior distribution and its computation, but asymptotic properties of the proportional hazard model have not been studied.

The asymptotic properties of the posterior are, however, an important theoretical issue in nonparametric Bayesian models, for there are many unexpected phenomena reported in the literature. For example, Diaconis and Freedman [4] showed that nonparametric posteriors could have inconsistency even with reasonable priors. They argued further that the inconsistency of the posterior in nonparametric problems is a rule, not an exception. The related work on this issue includes Ghosal, Ghosh and Ramamoorthi [8] and Barron, Schervish and Wasserman [1]. For right-censored data, Kim and Lee [13] showed that not all neutral to the right prior processes have consistent posteriors and gave sufficient conditions for consistency.

This unfortunate phenomenon continues to occur in the posterior convergence rate. See [2, 9, 23, 25]. These examples cast doubt on the Bernstein–von Mises theorem in nonparametric models, which states that the posterior distribution centered at the maximum likelihood estimator is asymptotically equivalent to the sampling distribution of the maximum likelihood estimator. See also [7]. In contrast, however, Shen [22] proved that even in semiparametric/nonparametric models, if the parameter of interest is of finite dimension, one does not need to worry because the Bernstein–von Mises theorem holds for finite-dimensional parameters.

If the Bernstein–von Mises theorem does not hold, it often implies the Bayesian credible set has zero efficiency relative to the frequentist confidence interval. The validity of the Bernstein–von Mises theorem also has an important implication in practice, because the Bernstein–von Mises theorem warrants use of Bayesian credible sets as frequentist confidence intervals asymptotically. Kim and Lee [16] studied the Bernstein–von Mises theorem for right-censored survival data without covariates. They found that for any $0 < \alpha \leq 1/2$ there is a consistent prior process neutral to the right whose posterior convergence rate is exactly $n^{-\alpha}$ and also showed that for popular prior processes such as beta, gamma and Dirichlet processes, indeed the Bernstein–von Mises theorem does hold.

In this paper we prove the Bernstein–von Mises theorem for Bayesian analysis of the proportional hazard model. The proof consists of the two Bernstein–von Mises theorems: one for the marginal posterior distribution of the regression coefficients and the other for the conditional posterior distribution of the baseline cumulative hazard functions given the regression coefficients. These two Bernstein–von Mises theorems together yield the Bernstein–von Mises theorem of the joint posterior distribution of the regression coefficients and the baseline cumulative hazard function. The main idea of the proof of the Bernstein–von Mises theorem of the marginal posterior distribution of the regression coefficients is to show that (i) on $1/\sqrt{n}$ neighborhoods of the true value, the posterior density converges to the targeted normal density with respect to the $L_1$ norm and (ii) outside $1/\sqrt{n}$ neighborhoods of the true value, the posterior mass vanishes eventually. For (i), we approximate the posterior distribution with the product of the partial likelihood and prior, and show that the product of the partial likelihood and prior converges to the target normal distribution with respect to the $L_1$ norm. The proof of (ii) is the harder part since the posterior distribution is not log-concave. For (ii), we use a sequence of log-concave functions which dominate the posterior distribution and whose total masses vanish eventually outside $1/\sqrt{n}$ neighborhoods of the true value. The proof of the Bernstein–von Mises theorem for the baseline cumulative hazard function given the regression coefficients exploits the functional central limit theorem for independent increment (II) processes (Theorem 19 of Section V.4 in [20]), for the conditional posterior distribution of the baseline cumulative hazard function given the regression coefficients is an II process.

The paper is organized as follows. In Section 2 prior processes neutral to the right are reviewed briefly and the posterior distribution of the regression coefficients and the baseline hazard function is given. In Section 3 the main results are stated and examples are given. Section 4 proves the main results with key lemmas, whose proofs are presented in the Appendix.

**2. Neutral to the right processes as priors.** The postulation of the proportional hazard model is as follows. Let $X_1, \ldots, X_n$ be survival times with covariates $Z_1, \ldots, Z_n$, where $Z_i \in R^p$, $i = 1, \ldots, n$. Suppose the distribution $F_i$ of $X_i$ with covariate $Z_i$ is given by

$$1 - F_i(t) = (1 - F(t))^{\exp(\beta^T Z_i)}$$

for an unknown regression parameter $\beta \in R^p$ and where $F$ is an unknown distribution of a survival time with covariate being 0. In most applications, the survival times are subject to right censoring, that is, $(T_1, \delta_1, Z_1), \ldots, (T_n, \delta_n, Z_n)$ are observed, where $T_i = \min(C_i, T_i)$, $\delta_i = I(X_i \leq C_i)$ and $C_1, \ldots, C_n$ are independent random variables with the common distribution function $G$.

In the proportional hazard model, there are two parameters: the regression coefficients $\beta$ and the baseline distribution function $F$. For prior distributions, we take a process neutral to the right [5] for $F$ and a usual parametric prior distribution for $\beta$. Processes neutral to the right include many popular prior processes such as Dirichlet processes, gamma processes and beta processes.

We say that a prior process on the c.d.f. $F$ is a process neutral to the right if the corresponding c.h.f. $A$ is a nondecreasing independent increment (NII) process such that $A(0) = 0$, $0 \leq \Delta A(t) \leq 1$ for all $t$ with probability 1 and either $\Delta A(t) = 1$ for some $t > 0$ or $\lim_{t \to \infty} A(t) = \infty$ with probability 1. See [5] for the original definition of processes neutral to the right and see [10, 12, 13] for the connection between the definition given here and Doksum's definition. From what follows, the term *NII process* is used for a prior process of the c.h.f. $A$ which induces a process neutral to the right on $F$.

The Lévy measure $\nu$ of an NII process $A$ is defined by

$$\nu([0,t] \times B) = \mathrm{E}\bigg(\sum_{s \in [0,t]} I(\Delta A(s) \in B \setminus \{0\})\bigg)$$

where $t \geq 0$ and $B$ is a Borel subset of $[0,1]$. Conversely, for any $\sigma$-finite measure $\nu$ defined on $[0,\infty) \times [0,1]$ which satisfies, for all $t > 0$, $\int_0^t \int_0^1 x\nu(ds,dx) < \infty$, there exists a unique NII process whose Lévy measure is $\nu$. Hence, any NII process can be characterized by its Lévy measure.

The mean and variance of an NII process $A(t)$ with Lévy measure $\nu$ can be conveniently calculated by the formulas

$$\mathrm{E}(A(t)) = \int_0^t \int_0^1 x\nu(ds,dx) \tag{1}$$

and

$$\mathrm{Var}(A(t)) = \int_0^t \int_0^1 x^2 \nu(ds,dx) - \sum_{s \leq t}\bigg(\int_0^1 x\nu(\{s\},dx)\bigg)^2. \tag{2}$$

These formulas constitute basic facts for the asymptotic theory of the posterior and will be used subsequently in this paper.

Let $q_n$ be the number of distinct uncensored observations and let $t_1 < t_2 < \cdots < t_{q_n}$ be the ordered distinct uncensored observations. Define two sets $D_n(t)$ and $R_n(t)$ by

$$D_n(t) = \{i : T_i = t, \delta_i = 1, i = 1, \ldots, n\}$$

and

$$R_n(t) = \{i : T_i \geq t, i = 1, \ldots, n\}.$$

Let $R_n^+(t) = R_n(t) - D_n(t)$.

A priori, let the baseline c.d.f. $F$ be a process neutral to the right such that the corresponding c.h.f. $A$ is an NII process with a Lévy measure $\nu$ of the form

$$\nu(dt, dx) = f_t(x)\, dx\, dt \tag{3}$$

for $x \in [0,1]$, and let $\pi(\beta)$ be the prior density function for $\beta$. Without loss of generality, we assume that $T_1 \leq T_1 \leq \cdots \leq T_n$. The next theorem provides the posterior distribution of $\beta$ as well as $A$. The proof is in [14].

THEOREM 2.1. *Let $D_n = ((T_1, \delta_1, Z_1), \ldots, (T_n, \delta_n, Z_n))$.*

(i) *Conditional on $\beta$ and $D_n$, the posterior distribution of $F$ is a process neutral to the right with Lévy measure*

$$\nu(dt, dx|\beta, D_n) = (1-x)^{\sum_{j \in R_n(t)} \exp(\beta^T Z_j)} f_t(x) \, dx \, dt \tag{4}$$

$$+ \sum_{i=1}^{q_n} dH_{ni}(x|\beta) \delta_{t_i}(dt),$$

*where $\delta_a$ is the point measure at $a$ and $H_{ni}(\cdot|\beta)$ is the probability measure defined on $[0,1]$ with density*

$$h_{ni}(x|\beta) \propto \left[ \prod_{j \in D_n(t_i)} (1 - (1-x)^{\exp(\beta^T Z_j)}) \right] \tag{5}$$

$$\times (1-x)^{\sum_{j \in R_n^+(t_i)} \exp(\beta^T Z_j)} f_{t_i}(x).$$

(ii) *The marginal posterior distribution of $\beta$ is*

$$\pi(\beta|D_n) \propto e^{-\rho_n(\beta)} \prod_{i=1}^{q_n} \int_0^1 h_{ni}(x|\beta) \, dx \, \pi(\beta), \tag{6}$$

*where*

$$\rho_n(\beta) = \sum_{i=1}^n \int_0^{T_i} \int_0^1 (1 - (1-x)^{\exp(\beta^T Z_i)})(1-x)^{\sum_{j=i+1}^n \exp(\beta^T Z_j)} f_t(x) \, dx \, dt,$$

*for $j = 1, \ldots, n$ and $\sum_{j=i+1}^n \exp(\beta^T Z_j) = 0$ when $i = n$.*

**3. Main result.** Let $\beta_0$ and $F_0$ be the true values of the parameters where $X_1, \ldots, X_n$ are generated and let $A_0$ be the cumulative hazard function of $F_0$. In this section we present the Bernstein–von Mises theorem of the posterior distribution of $(\beta, A)$. That is, we show that asymptotically the posterior distribution of $(\beta, A)$ centered at the maximum partial likelihood estimator (MLE) $(\hat{\beta}, \hat{A})$ is the same as the asymptotic distribution of the MLE itself.

The following conditions are assumed to hold in the remainder of this section:

(A1) $A_0$ is absolutely continuous.
(A2) For a positive constant $\tau$, $F_0(\tau) < 1$, $G(\tau-) < 1$ and $G(\tau) = 1$.
(A3) $Z_1, \ldots, Z_n$ are i.i.d. $p$-dimensional random vectors such that $\|Z_1\| \leq M_z < \infty$ with probability 1 for some constant $M_z$ where

$$\|Z_1\| = \sum_{i=1}^p |Z_{1i}|.$$

(A4) If $\Pr(c'Z_1 = 0) = 1$, then $c = 0$.
(A5) $\pi(\beta)$ is continuous at $\beta_0$ with $\pi(\beta_0) > 0$.

Condition (A1) prevents ties. If $A_0$ has a finite number of discontinuity points, then the proof can be done separately on the continuous part and the discrete part. Condition (A2) assumes that some patients remain in the study until time $\tau$, which is necessary to recover the information of $A(t)$ on $[0, \tau]$. If condition (A2) holds for all $\tau$, the Bernstein–von Mises theorem for $A$ holds on $[0, \infty)$. But, note that even if $\tau < \infty$, the Bernstein–von Mises theorem for $\beta$ holds as long as $\tau > 0$. Condition (A3) is for technical purposes, and condition (A4) is to avoid collinearity among the covariates. Condition (A5) is a standard assumption for Bernstein–von Mises type results.

Let $\hat{\beta}$ be the maximum (partial) likelihood estimator which maximizes the partial likelihood

$$L_n(\beta) = \prod_{i=1}^{q_n} \prod_{j \in D_n(t_i)} \frac{\exp(\beta^T Z_j)}{n^{-1} \sum_{k \in R_n(t_i)} \exp(\beta^T Z_k)}.$$

Let

$$\hat{A}(t) = \int_0^t \frac{dN.(s)}{\sum_{i \in R_n(s)} \exp(\hat{\beta}^T Z_i)},$$

where $N.(t) = \sum_{i=1}^n N_i(t)$ and $N_i(t) = I(T_i \leq t, \delta_i = 1)$. In fact, $\hat{A}$ is Breslow's estimator of the baseline hazard function [3]. We introduce the notation

$$U_0(t) = \int_0^t \frac{dA_0(s)}{S^0(s : \beta_0)},$$

$$e_0(t) = \int_0^t \frac{S^1(s : \beta_0) \, dA_0(s)}{S^0(s : \beta_0)},$$

$$S^0(t : \beta) = \mathrm{E}(\exp(\beta^T Z_1) I(T_1 \geq t)),$$

$$S^1(t : \beta) = \mathrm{E}(Z_1 \exp(\beta^T Z_1) I(T_1 \geq t)),$$

$$I(\beta) = \int_0^\tau V(t : \beta) S^0(t : \beta) \, dA_0(t),$$

$$V(t : \beta) = S^2(t : \beta) / S^0(t : \beta) - e_0(t)^2,$$

$$S^2(t : \beta) = \mathrm{E}(Z_1 Z_1^T \exp(\beta^T Z_1) I(T_1 \geq t)).$$

Assume that a priori $A$ is an NII process with Lévy measure given by

$$(7) \qquad \nu(ds, dx) = \frac{g_s(x)}{x} \, dx \, \lambda(s) \, ds, \qquad s \geq 0, 0 \leq x \leq 1,$$

where $\int_0^1 g_t(x) \, dx = 1$ for all $t \in [0, \tau]$ and that $\lambda(t)$ is bounded and positive on $(0, \tau)$.

REMARK. Comparing (3) and (7), we can see that

$$\int_0^t \lambda(s)\,ds = \int_0^t \int_0^1 x f_s(x)\,dx\,ds = \mathrm{E}(A(t))$$

and $g_t(x) = x f_t(x)/\lambda(t)$ provided $\lambda(t) > 0$.

REMARK. Positiveness of $\lambda(t)$ on $t \in (0, \tau)$ is necessary for the Bernstein–von Mises theorem. Suppose $\lambda(t) = 0$ for $t \in [c,d]$ where $0 < c < d < \tau$. Then both the prior and posterior put mass 1 on the set of c.h.f.'s, $A$ with $A(d) = A(c)$.

For the Bernstein–von Mises theorem, we need the following two conditions:

(C1) There exists a positive number $\varsigma$ such that

$$\sup_{t \in [0,\tau], x \in [0,1]} (1-x)^{1-\varsigma} g_t(x) < \infty.$$

(C2) There exists a function $k(t)$ defined on $[0, \tau]$ such that for some $\alpha > 1/2$ and $\varepsilon > 0$

$$\sup_{t \in [0,\tau], h \in [0,\varepsilon]} \left| \frac{g_t(h) - k(t)}{h^\alpha} \right| < \infty$$

and

$$0 < \inf_{t \in [0,\tau]} k(t) \le \sup_{t \in [0,\tau]} k(t) < \infty.$$

Throughout this paper, we let

$$g^* = \sup_{t \in [0,\tau], x \in [0,1]} (1-x)^{1-\varsigma} g_t(x),$$

$k_* = \inf_{t \in [0,\tau]} k(t)$ and $k^* = \sup_{t \in [0,\tau]} k(t)$.

Conditions (C1) with $\varsigma = 0$ and (C2) are used for the Bernstein–von Mises theorem of the survival function without covariates by Kim and Lee [16]. The most delicate part of the proof in this paper is to show that the tail probability of the posterior distribution of $\beta$ converges to 0 sufficiently fast. The positiveness of $\varsigma$ plays an important role for this.

The following theorems are the main results of this paper. We first state the result, an interesting result in its own right, that the marginal posterior density of $\beta$ converges to a normal density in the $L_1$ norm. This is stronger than the usual Bernstein–von Mises theorem, which states that the posterior converges weakly to a normal distribution in probability, because our result states that the posterior density converges to a normal density in the $L_1$ norm with probability 1.

THEOREM 3.1. *Under conditions* (C1) *and* (C2),

$$\lim_{n \to \infty} \|f_n - \phi\| = 0 \tag{8}$$

with probability 1, where $f_n$ is the posterior density of $\sqrt{n}(\beta - \hat{\beta})$, $\phi$ is the normal density with mean 0 and variance $I(\beta_0)^{-1}$, and $\|\cdot\|$ is the $L_1$ norm.

The next theorem states that the conditional distribution of $\sqrt{n}(A - \hat{A})$ given $\beta$ and data converges to a Gaussian process.

THEOREM 3.2. *Under conditions* (C1) *and* (C2),

$$\mathcal{L}(\sqrt{n}(A(\cdot) - \hat{A}(\cdot))|\sqrt{n}(\beta - \hat{\beta}) = x, D_n) \xrightarrow{d} W(U_0(\cdot)) - xe_0(\cdot)$$

on $D[0,\tau]$ with probability 1, where $W$ is standard Brownian motion. Here, $D[0,\tau]$ is the space of right-continuous functions on $[0,\tau]$ with left limits existing on $[0,\tau]$ equipped with the uniform topology.

The proofs of Theorems 3.1 and 3.2 are presented in Section 4. Combining Theorems 3.1 and 3.2, we can prove the main theorem stated below.

THEOREM 3.3. *Under conditions* (C1) *and* (C2),

$$\mathcal{L}(\sqrt{n}(A(\cdot) - \hat{A}(\cdot), \beta - \hat{\beta})|D_n) \xrightarrow{d} (W(U_0(\cdot)) - Xe_0(\cdot), X) \tag{9}$$

as $n \to \infty$ on $D[0,\tau] \times R^p$ where $X$ is a multivariate normal random vector with mean 0 and variance $I^{-1}(\beta_0)$ and $W$ is standard Brownian motion independent of $X$.

PROOF. Theorems 3.1 and 3.2 prove the convergence of the marginal posterior distribution of $\beta$ and the conditional posterior distribution of $A$ given $\beta$. To prove the convergence of the joint posterior distribution of $\beta$ and $A$, note that Theorem 3.1 implies the strong convergence of the marginal posterior distribution of $\sqrt{n}(\beta - \hat{\beta})$ to the distribution of $X$. Applying Theorem 2 of [21], we complete the proof. □

REMARK. It should be noted that the limiting distribution (9) is the same as that of the maximum likelihood estimators centered at the true values.

REMARK. From (9), we can see that marginally the posterior distributions of $\sqrt{n}(\beta - \hat{\beta})$ and $\sqrt{n}(A - \hat{A})$ converge weakly to a normal distribution and a Gaussian process, respectively.

In the following examples, we show that most popular prior processes such as beta processes and gamma processes satisfy condition (C1). For condition (C2), see [16].

EXAMPLE 1 (Beta process). The beta process with mean $\Lambda$ and scale parameter $c$ is an NII process with Lévy measure $\nu$ given by

$$\nu(dt, dx) = \frac{c(t)}{x}(1-x)^{c(t)-1}\, dx\, d\Lambda(t).$$

Suppose that $\Lambda(t)$ is absolutely continuous with $\lambda(t) = d\Lambda(t)/dt$. Then $g_t(x) = c(t)(1-x)^{c(t)-1}$. If $\inf_{t \in [0,\tau]} c(t) > 0$ and $\sup_{t \in [0,\tau]} c(t) < \infty$, condition (C1) holds with $\varsigma = \inf_{t \in [0,\tau]} c(t)$.

EXAMPLE 2 (Gamma process). A priori, assume that $Y(t) = -\log(1 - F(t))$ is a gamma process with parameters $(\Lambda(t), c(t))$ with $\Lambda(t) = \int_0^t \lambda(s)\, dx$, where $\lambda(t)$ is a positive bounded function on $t \in (0, \tau)$. Furthermore, assume that $c(t)$ is continuous around $t = 0$ and $0 < \inf_{t \in [0,\tau]} c(t)(= c_*) \leq \sup_{t \in [0,\tau]} c(t)(= c^*) < \infty$. Here, the gamma process with parameters $(\Lambda(t), c(t))$ is defined by

$$Y(t) = \int_0^t \frac{1}{c(s)}\, dX(s),$$

where $X(t)$ is an NII process whose marginal distribution of $X(t)$ is a gamma distribution with parameters $(\int_0^t c(s)\, d\Lambda(s), 1)$. For details of this definition, see [19]. This prior process was used by Doksum [5], Ferguson and Phadia [6] and Kalbfleisch [11]. Since

$$\log \mathrm{E}(\exp(-\theta Y(t))) = \int_0^t \int_0^\infty (e^{-\theta x} - 1)\frac{c(s)}{x}\exp(-c(s)x)\, dx\, d\Lambda(s),$$

it can be shown that the c.h.f. $A$ of $F$ is an NII process with Lévy measure given by

$$\nu(ds, dx) = \tilde{c}(s)\frac{1}{-\log(1-x)}(1-x)^{c(s)-1}\, dx\, d\tilde{\Lambda}(s),$$

where

$$\tilde{c}(t) = \left(\int_0^1 \frac{x}{-\log(1-x)}(1-x)^{c(t)-1}\, dx\right)^{-1}$$

and

$$\tilde{\Lambda}(t) = \int_0^t \frac{c(s)}{\tilde{c}(s)}\, d\Lambda(s).$$

Therefore, we have

$$g_t(x) = \tilde{c}(t)\frac{x}{-\log(1-x)}(1-x)^{c(t)-1}, \qquad 0 \leq x \leq 1.$$

Now,

$$g_t(x) = \tilde{c}(t)\frac{x(1-x)^{c_*/2}}{-\log(1-x)}(1-x)^{c(t)-c_*/2-1}$$

$$\leq \left(\sup_{t\in[0,\tau]} \tilde{c}(t)\right) m(1-x)^{c(t)-c_*/2-1},$$

where

$$m = \sup_{t\in[0,\tau]} \frac{x(1-x)^{c_*/2}}{-\log(1-x)}.$$

It is easy to show that $\sup_{t\in[0,\tau]} \tilde{c}(t) < \infty$ and so condition (C1) follows with $\varsigma < c_*/2$.

**4. Proof of the main results.** For a given sequence of random variables $Z_n$, we write $Z_n = O(n^\delta)$ with probability 1 if there exists a constant $M > 0$ such that $Z_n/n^\delta \leq M$ for all but finitely many $n$ with probability 1. Also, we write $Z_n = o(n^\delta)$ with probability 1 if $Z_n/n^\delta$ converges to 0 with probability 1. For a given finite-dimensional array of real numbers $C$, $\|C\|$ is defined as the sum of all the absolute values of the elements of $C$.

Let $d(i)$ be the integer such that $T_{d(i)} = t_i$ and $\delta_{d(i)} = 1$. Note that since we assume that the true distribution $F_0$ is continuous there is no tie among the uncensored observations and so $d(i)$ is well defined.

4.1. *Proof of Theorem* 3.1. Let

$$h_n(\beta) = -\rho_n(\beta) + \sum_{i=1}^{q_n} \log\left(n\int_0^1 h_{ni}(x|\beta)\, dx\right).$$

Then we have

$$\pi(\beta|D_n) \propto \exp(h_n(\beta))\pi(\beta),$$

and so the posterior density of $\sqrt{n}(\beta - \hat{\beta})$ becomes $f_n(h) = g_n(h)/C_n$ where

$$g_n(h) = \exp(h_n(\hat{\beta} + h/\sqrt{n}) - h_n(\hat{\beta}))\pi(\hat{\beta} + h/\sqrt{n})$$

and $C_n = \int_{R^p} g_n(h)\, dh$. Hence, the proof of Theorem 3.1 will be completed if we prove that

$$(10) \qquad \int_{R^p} |g_n(h) - \psi(h)\pi(\beta_0)|\, dh \to 0$$

with probability 1, where

$$\psi(h) = \exp(-h^T I(\beta_0) h/2).$$

Let $l_n(\beta) = \log L_n(\beta)$ and define

$$\tilde{l}(\beta) = \beta^T \mathrm{E}(Z_1 I(\delta_1 = 1))$$
$$- \int_0^\tau \log(\mathrm{E}(e^{\beta^T Z_1} I(T_1 \geq t))) \mathrm{E}(e^{\beta_0^T Z_1} I(T_1 \geq t)) \, dA_0(t).$$

It is not hard to see that $\tilde{l}(\beta)$ is strictly concave with attainment of its maximum at $\beta_0$. Hence,

$$\sup_{\beta \in B} |l_n(\beta)/n - \tilde{l}(\beta)| \to 0 \tag{11}$$

with probability 1 for any compact subset $B$ of $R^p$.

We recall the following properties of $\hat{\beta}$ and $l_n(\beta)$ from [24] or [15]. First, $\hat{\beta}$ is consistent (i.e., $\hat{\beta} \to \beta_0$ with probability 1). Let $l_n^{(k)}(\beta)$ be the $k$th derivative of $l_n(\beta)$ in $\beta$. Then $-l_n^{(2)}(\hat{\beta})/n \to I(\beta_0)$ with probability 1. Also, $\sup_{\beta \in B} \|l_n^{(3)}(\beta)\| = O(n)$ with probability 1 for any compact subset $B$ of $R^p$.

We need the following two lemmas whose proofs are in the Appendix.

LEMMA 1. *For any compact subset $B$ of $R^p$,*

$$\sup_{\beta \in B} \frac{1}{n} \|h_n^{(k)}(\beta) - l_n^{(k)}(\beta)\| = o(1)$$

*with probability 1 for $k = 0, 1, 2, 3$.*

LEMMA 2.

$$\|h_n^{(1)}(\hat{\beta})\| = o(\sqrt{n})$$

*with probability 1.*

We decompose (10) by

$$\int_{R^p} |g_n(h) - \psi(h)\pi(\beta_0)| \, dh \leq \int_{|h| \leq K} |g_n(h) - \psi(h)\pi(\beta_0)| \, dh \tag{12}$$

$$+ \int_{|h| > K} \psi(h)\pi(\beta_0) \, dh \tag{13}$$

$$+ \int_{K < |h| \leq \sqrt{n}\delta} g_n(h) \, dh \tag{14}$$

$$+ \int_{|h| > \sqrt{n}\delta} g_n(h) \, dh. \tag{15}$$

We will show that for given $\varepsilon > 0$, there exist positive constants $K$ and $\delta$ such that the four terms become smaller than $\varepsilon$ for all sufficiently large $n$.

For (12), we will exploit the standard techniques used for the proof of the Bernstein–von Mises theorem for parametric models. First, using Taylor expansion, we write

$$\log(g_n(h)) = \frac{h^T}{\sqrt{n}} h_n^{(1)}(\hat{\beta}) - \frac{1}{2} h^T \left( -\frac{1}{n} h_n^{(2)}(\hat{\beta}) \right) h + R_n(h) \tag{16}$$
$$+ \log(\pi(\hat{\beta} + h/\sqrt{n})),$$

where $h_n^{(k)}$ is the $k$th derivative of $h_n$ in $\beta$. Lemma 2 implies, for all $h$,

$$\frac{h}{\sqrt{n}} h_n^{(1)}(\hat{\beta}) \to 0. \tag{17}$$

Lemma 1 with the properties of $l_n(\beta)$ yields, for all $h$,

$$-\frac{1}{2} h^T \left( -\frac{1}{n} h_n^{(2)}(\hat{\beta}) \right) h \to -\frac{1}{2} h^T I(\beta_0) h, \tag{18}$$

$$|R_n(h)| \to 0. \tag{19}$$

Also, we have

$$\log(\pi(\hat{\beta} + h/\sqrt{n})) \to \pi(\beta_0) \tag{20}$$

uniformly on $\{|h| \leq K\}$ with probability 1 for any $K > 0$. Now,

$$|g_n(h) - \psi(h)\pi(\beta_0)|$$
$$\leq |g_n(h) - \psi(h)\pi(\hat{\beta} + h/\sqrt{n})| + |\psi(h)\pi(\hat{\beta} + h/\sqrt{n}) - \psi(h)\pi(\beta_0)|$$
$$\leq \psi(h)\pi(\hat{\beta} + h/\sqrt{n})$$
$$\times \left| \exp\left( h^T \frac{h_n^{(1)}(\hat{\beta})}{\sqrt{n}} - \frac{1}{2} h^T \left( -\frac{1}{n} h_n^{(2)}(\hat{\beta}) - I(\beta_0) \right) h + R_n(h) \right) - 1 \right|$$
$$+ \psi(h)\pi(\beta_0) \left| \frac{\pi(\hat{\beta} + h/\sqrt{n})}{\pi(\beta_0)} - 1 \right|.$$

By (17)–(20), we get, for any $K > 0$,

$$\sup_{|h| \leq K} |g_n(h) - \psi(h)\pi(\beta_0)| \to 0$$

with probability 1, and thus

$$\int_{|h| \leq K} |g_n(h) - \psi(h)\pi(\beta_0)| \, dh \to 0$$

with probability 1.

We can make (13) as small as possible by choosing $K$ sufficiently large.

As for (14), note that by Lemma 1 with the property of $l_n(\beta)$, there exists a constant $M$ such that $\sup_{|\beta|\leq K}\|h_n^{(3)}(\beta)/n\| \leq M$ for sufficiently large $n$. Hence, we can write

$$R_n(h) \leq \sum_{i,j,k=1}^{p} \frac{|h_i||h_j||h_k|}{6\sqrt{n}} \left(\left\|\frac{1}{n}h_n^{(3)}(\tilde{\beta})\right\|\right)$$

(21)

$$\leq \frac{p^2\delta}{6} M h^T h,$$

for some $\tilde{\beta}$ in between $\beta_0$ and $\hat{\beta}$. Let $\eta > 0$ be the smallest eigenvalue of $I(\beta_0)$. Since $-h_n^{(2)}(\hat{\beta})/n \to I(\beta_0)$ with probability 1, we have

(22) $$h^T(-h_n^{(2)}(\hat{\beta})/n)h \geq (\eta - o(1))h^T h.$$

Also, when $|h| > 1$,

(23) $$h^T \frac{h_n^{(1)}(\hat{\beta})}{\sqrt{n}} \leq |h|\left\|\frac{h_n^{(1)}(\hat{\beta})}{\sqrt{n}}\right\| = |h|^2\frac{o(1)}{|h|} \leq o(1)h^T h$$

with probability 1. Now, combining (21), (22) and (23), we have

$$h_n(\hat{\beta} + h/\sqrt{n}) - h_n(\hat{\beta}) \leq -h^T h(\eta/2 - p^2\delta M/6 + o(1))$$

when $|h| > 1$ for all sufficiently large $n$ with probability 1. Set $\delta$ sufficiently small that $\eta/2 - p^2\delta M/6 (= \kappa) > 0$ and $\sup_{|\beta-\beta_0|\leq 2\delta}\pi(\beta)(= \varrho) < \infty$. Then

$$\int_{K\leq |h|\sqrt{n}\delta} g_n(h)\,dh \leq \int_{K\leq |h|\leq \sqrt{n}\delta} \exp(-|h|^2(\kappa + o(1)))\pi(\hat{\beta} + h/\sqrt{n})\,dh$$

$$\leq \varrho \int_{K\leq |h|\leq \sqrt{n}\delta} \exp(-|h|^2\kappa/2)\,dh$$

for all sufficiently large $n$ with probability 1. Hence, we can make (14) as small as possible by choosing a sufficiently large $K$.

For (15), let $\psi(x) = \int_0^1 (1-(1-y)^{x-1})/y\,dy$. Then it can be shown that $\sup_{\varsigma\leq x<\infty} x\psi'(x) = \psi^* < \infty$ where $\psi'(x) = d\psi(x)/dx$. Note that

$$h_n(\beta) \leq \sum_{i=1}^{q_n} \log\biggl(n\int_0^1 \frac{1-(1-x)^{\exp(\beta^T Z_{d(i)})}}{x}$$

$$\times (1-x)^{\sum_{j\in R_n^+(t_i)}\exp(\beta^T Z_j)} g_{t_i}(x)\,dx\biggr)$$

$$\leq \sum_{i=1}^{q_n} \log\biggl(g^* n\int_0^1 \frac{1-(1-x)^{\exp(\beta^T Z_{d(i)})}}{x}$$

$$\times (1-x)^{\sum_{j\in R_n^+(t_i)} \exp(\beta^T Z_j)+\varsigma -1} dx\bigg)$$

$$\leq \sum_{i=1}^{q_n} \log\bigg[g^* n\bigg(\psi\bigg(\sum_{j\in R_n(t_i)} \exp(\beta^T Z_j)+\varsigma\bigg) -\psi\bigg(\sum_{j\in R_n^+(t_i)} \exp(\beta^T Z_j)+\varsigma\bigg)\bigg)\bigg]$$

$$(24) \qquad \leq \sum_{i=1}^{q_n} \log\bigg(g^*\psi^* n \frac{\exp(\beta^T Z_{d(i)})}{\sum_{j\in R_n^+(t_i)} \exp(\beta^T Z_j)+\varsigma}\bigg)$$

$$(25) \qquad \leq Cq_n + \sum_{i=1}^{q_n} \log\bigg(n \frac{\exp(\beta^T Z_{d(i)})}{\sum_{j\in R_n^+(t_i)} \exp(\beta^T Z_j)}\bigg),$$

where $C = \log(g^*\psi^*)$. Here the inequality in (24) follows from

$$\psi\bigg(\sum_{j\in R_n(t_i)} \exp(\beta^T Z_j)+\varsigma\bigg) -\psi\bigg(\sum_{j\in R_n^+(t_i)} \exp(\beta^T Z_j)+\varsigma\bigg)$$

$$= \exp(\beta^T Z_{d(i)})\psi'(a)$$

$$= \frac{\exp(\beta^T Z_{d(i)})}{a} a\psi'(a)$$

$$\leq \frac{\exp(\beta^T Z_{d(i)})}{\sum_{j\in R_n^+(t_i)} \exp(\beta^T Z_j)+\varsigma}\psi^*,$$

where $a$ is a positive number between $\sum_{j\in R_n^+(t_i)} \exp(\beta^T Z_j) + \varsigma$ and $\sum_{j\in R_n(t_i)} \exp(\beta^T Z_j) + \varsigma$.

Let

$$l_n^+(\beta) = \sum_{i=1}^{q_n} \log\bigg(n\frac{\exp(\beta^T Z_{d(i)})}{\sum_{j\in R_n^+(t_i)} \exp(\beta^T Z_j)}\bigg).$$

Note that $R_n^+(t_i)$ are nonempty sets and so $l_n^+(\beta)$ is well defined for all sufficiently large $n$ with probability 1. Also, by direct calculation we can see that $l_n^+(\beta)$ is a strictly concave function. Since $\sup_{\beta\in B} |l_n(\beta) - l_n^+(\beta)| = O(1)$ for any compact subset $B$ of $R^p$, we have $\sup_{\beta\in B} |l_n^+(\beta)/n - \tilde{l}(\beta)| \to 0$. Now, choose $m$ such that

$$\sup_{\beta:|\beta-\beta_0|=m} (\tilde{l}(\beta) - \tilde{l}(\beta_0)) \leq -qC - \eta$$

for some $\eta > 0$ where $q = \Pr\{T_1 \leq \tau, \delta_1 = 1\}$. Then
$$\sup_{\beta:|\beta-\beta_0|\geq m} \frac{l_n^+(\beta)}{n} \leq \sup_{\beta:|\beta-\beta_0|=m} \frac{l_n^+(\beta)}{n} = \sup_{\beta:|\beta-\beta_0|=m} \tilde{l}(\beta) + o(1).$$

Since $q_n/n \to q$ and $h_n(\hat{\beta})/n \to \tilde{l}(\beta_0)$ with probability 1 by Lemma 1, we have
$$\sup_{\beta:|\beta-\beta_0|\geq \delta/2} l_n^+(\beta)/n - h_n(\hat{\beta})/n$$
$$\leq \sup_{\beta:\delta/2 \leq |\beta-\beta_0|\leq m} (l_n^+(\beta)/n - h_n(\hat{\beta})/n)$$
$$+ \sup_{\beta:|\beta-\beta_0|\geq m} (l_n^+(\beta)/n - h_n(\hat{\beta})/n)$$
$$= -qC - \eta + o(1).$$

Finally,
$$\int_{|h|\geq \sqrt{n}\delta} g_n(h)\,dh = n^{p/2} \int_{|\beta-\hat{\beta}|\geq \delta} e^{h_n(\beta)-h_n(\hat{\beta})} \pi(\beta)\,d\beta$$
$$\leq n^{p/2} \sup_{\beta:|\beta-\hat{\beta}|\geq \delta} e^{h_n(\beta)-h_n(\hat{\beta})}$$
$$= n^{p/2} \exp\left[n\left(q_n C/n + \sup_{\beta:|\beta-\beta_0|\geq \delta/2} l_n^+(\beta)/n - h_n(\hat{\beta})/n\right)\right]$$
$$\leq n^{p/2} \exp[n(-\eta + o(1))]$$
$$\leq n^{p/2} e^{-n\eta/2} \to 0$$

for all sufficiently large $n$ with probability 1 and the proof is done.

4.2. *Proof of Theorem 3.2.* Let $\theta_n = (\sqrt{n}(\beta - \hat{\beta}) = h, D_n)$ be given. We decompose $\sqrt{n}(A(\cdot) - \hat{A}(\cdot))$ by

(26) $$\sqrt{n}(A(\cdot) - \hat{A}(\cdot)) = \sqrt{n}(A(\cdot) - A^J(\cdot))$$
(27) $$+ \sqrt{n}(A^J(\cdot) - \tilde{A}(\cdot))$$
(28) $$+ \sqrt{n}(\tilde{A}(\cdot) - \hat{A}_h(\cdot))$$
(29) $$+ \sqrt{n}(\hat{A}_h(\cdot) - \hat{A}(\cdot)),$$

where
$$A^J(t) = \sum_{i=1}^{q_n} \Delta A(t_i) I(t_i \leq t),$$
$$\tilde{A}(t) = \mathrm{E}(A^J(t)|\theta_n)$$

and

$$\hat{A}_h(t) = \int_0^t \frac{dN.(u)}{\sum_{j \in R_n(u)} \exp(\hat{\beta}_h^T Z_j)}$$

with $\hat{\beta}_h = \hat{\beta} + h/\sqrt{n}$ and $N.(t) = \sum_{i=1}^n I(T_i \leq t, \delta_i = 1)$. Then we will prove that when $\theta_n$ is given, with probability 1, (26) and (28) converge to 0 weakly, (27) converges to $W(U_0(\cdot))$ weakly, and (29) converges to $he_0(\cdot)$ on $D[0, \tau]$ with probability 1. Then Slutsky's theorem completes the proof.

For (26), let $\hat{\beta}_h = \hat{\beta} + h/\sqrt{n}$ and $Y_n(t) = A(t) - A^J(t)$. Then Theorem 2.1 yields that conditional on $\theta_n$, $Y_n(t)$ is an NII process with Lévy measure $\nu_{Y_n}$ given by

$$\nu_{Y_n}(dt, dx) = (1-x)^{\sum_{j \in R_n(t)} \exp(\hat{\beta}_h^T Z_j)} \frac{g_t(x)}{x} dx \, \lambda(t) \, dt.$$

Since $Y_n$ is nondecreasing, $\sup_{t \in [0, \tau]} |\sqrt{n} Y_n(t)| = \sqrt{n} Y_n(\tau)$ and so it suffices to show that $\mathcal{L}(\sqrt{n} Y_n(\tau) | \theta_n) \xrightarrow{d} 0$ with probability 1, which is equivalent to $\Pr\{|\sqrt{n} Y_n(\tau)| \geq \varepsilon | \theta_n\} \to 0$ with probability 1 for any $\varepsilon > 0$. By the Chebyshev inequality, we have

$$\Pr\{|\sqrt{n} Y_n(\tau)| \geq \varepsilon | \theta_n\} \leq \frac{1}{\varepsilon^2}((\sqrt{n} E(Y_n(\tau) | \theta_n))^2 + n \operatorname{Var}(Y_n(\tau) | \theta_n)).$$

Let

$$\phi_k = n^{k/2} \int_0^\tau \int_0^1 x^{k-1} (1-x)^{\sum_{j \in R_n(s)} \exp(\hat{\beta}_h^T Z_j)} g_s(x) \, dx \, \lambda(s) \, ds.$$

Then

$$\phi_k \leq n^{k/2} g^* \int_0^\tau \int_0^1 x^{k-1} (1-x)^{\sum_{j \in R_n(\tau)} \exp(\hat{\beta}_h^T Z_j) + \varsigma - 1} \, dx \, \lambda(s) \, ds$$

$$= n^{k/2} g^* \int_0^\tau \lambda(s) \, ds O(n^{-k}) = O(n^{-k/2}).$$

Since $\sqrt{n} E(Y_n(\tau) | \theta_n) = \phi_1$ and $n \operatorname{Var}(Y_n(\tau) | \theta_n) = \phi_2$, (26) converges to 0 on $D[0, \tau]$ with probability 1.

For (27)–(29), we need the following lemma, whose proof is in the Appendix.

LEMMA 3. *For any compact subset $B$ of $R^p$ and any positive integer $k$,*

$$\sup_{\beta \in B, 1 \leq i \leq q_n} \left| E(\Delta A^k(t_i) | \beta, D_n) - \frac{k! \Gamma(\sum_{j \in R(t_i)} \exp(\beta^T Z_j) + 1)}{\Gamma(\sum_{j \in R(t_i)} \exp(\beta^T Z_j) + k + 1)} \right|$$

$$= o(n^{-(k+1/2)})$$

*with probability* 1.

For (27), since conditional on $\theta_n$ the process $\sqrt{n}(A^J - \tilde{A})$ is an independent increment process, we utilize Theorem 19 of Section V.4 in [20]. Let $Y_n = \sqrt{n}(A^J - \tilde{A})$. We first prove the convergence of the finite-dimensional distribution by showing Lyapounov's condition. Suppose $0 \leq s < t \leq \tau$ are given. Note that

$$Y_n(t) - Y_n(s) = \sum_{s < t_i \leq t} \sqrt{n}(\Delta A(t_i) - \Delta \tilde{A}(t_i)).$$

Lemma 3 implies

$$\sup_{i=1,\ldots,q_n} \mathrm{E}[(\sqrt{n}(\Delta A(t_i) - \Delta \tilde{A}(t_i)))^4 | \theta_n] = O(n^{-2})$$

with probability 1. Hence

$$(30) \quad \sum_{s < t_i \leq t} \mathrm{E}[(\sqrt{n}(\Delta A(t_i) - \Delta \tilde{A}(t_i)))^4 | \theta_n] = \int_s^t O(n^{-2}) \, dN.(u) \to 0$$

with probability 1. Similarly, we have that

$$(31) \quad \begin{aligned} &\mathrm{Var}(Y_n(t) - Y_n(s) | \theta_n) \\ &= \sum_{s < t_i \leq t} \mathrm{E}\{[\sqrt{n}(\Delta A(t_i) - \Delta \tilde{A}(t_i))]^2 | \theta_n\} \\ &= \int_s^t \frac{n}{\sum_{j \in R_n(u)} \exp(\hat{\beta}_h^T Z_j)} (1 + o(n^{-1/2})) \frac{dN.(u)}{\sum_{j \in R_n(u)} \exp(\hat{\beta}_h^T Z_j)} \end{aligned}$$

with probability 1. Hence, Lemma A2 in [24] yields

$$(32) \quad \mathrm{Var}(Y_n(t) - Y_n(s) | \theta_n) \to U_0(t) - U_0(s)$$

with probability 1. Now (30) and (32) imply the finite-dimensional distributions of $Y_n$ converge to those of $W(U_0)$ weakly. Finally, note that

$$\Pr\{|Y_n(t) - Y_n(s)| \geq \varepsilon | \theta_n\} \leq \frac{1}{\varepsilon^2} \mathrm{Var}(Y_n(t) - Y_n(s) | \theta_n).$$

By (32), we have

$$\mathrm{Var}(Y_n(t) - Y_n(s) | \theta_n) = U_0(t) - U_0(s) + o(1)$$

with probability 1. Since $U_0(t)$ is continuous, with probability 1 we can make $\Pr\{|Y_n(t) - Y_n(s)| \geq \varepsilon | \theta_n\}$ as small as possible for all sufficiently large $n$ by choosing $t$ and $s$ sufficiently close. Hence Theorem 19 of Section V.4 in [20] allows us to conclude that $\mathcal{L}(Y_n | \theta_n)$ converges weakly to $W(U_0)$ on $D[0, \tau]$ with probability 1.

For (28), Lemma 3 yields that

$$\Delta \tilde{A}(t_i) = \frac{1}{\sum_{j \in R_n(t_i)} \exp(\hat{\beta}_h^T Z_j) + 1} + o(n^{-3/2}).$$

Therefore

$$\sup_{t\in[0,\tau]} |\sqrt{n}(\tilde{A}(t) - \hat{A}_h(t))|$$

$$= \int_0^\tau \sqrt{n} \left| \frac{\sum_{j\in R_n(u)} \exp(\hat{\beta}_h^T Z_j)}{\sum_{j\in R_n(u)} \exp(\hat{\beta}_h^T Z_j) + 1} \right.$$

$$\left. - 1 + \left(\sum_{j\in R_n(u)} \exp(\hat{\beta}_h^T Z_j)\right) o(n^{-3/2}) \right| \frac{dN.(u)}{\left|\sum_{j\in R_n(u)} \exp(\hat{\beta}_h^T Z_j)\right|}$$

$$= \int_0^\tau o(1) \frac{dN_n(u)}{\sum_{j\in R_n(u)} \exp(\hat{\beta}_h^T Z_j)} \to 0$$

with probability 1 by Lemma A.2 of [24].

Finally, the proof of (29) converging to 0 can be found in the proof of Theorem 3 in [15].

## APPENDIX: PROOF OF LEMMAS IN SECTION 4

LEMMA A.1. *Let*

$$\eta_i(x,\beta) = \frac{(1-(1-x)^{\exp(\beta^T Z_{d(i)})}) g_{t_i}(x)(1-x)}{k(t_i)x}$$

(33)

$$- \exp(\beta^T Z_{d(i)})(1-x)^{\exp(\beta^T Z_{d(i)})}$$

*and let $\eta_i^{(k)}(x,\beta)$ be the kth derivative of $\eta_i(x,\beta)$ in $\beta$. Let $\alpha' = \min\{1,\alpha\}$ where $\alpha$ is in condition* (C2). *Then for any compact subset $B$ of $R^p$, there exist constants $M_k$, $k=0,1,2,3$, such that*

$$\sup_{\beta\in B, x\in(0,1), 1\leq i\leq q_n} \left\| \frac{\eta_i^{(k)}(x,\beta)}{x^{\alpha'}} \right\| < M_k \tag{34}$$

*with probability* 1.

PROOF. Write

$$\eta_i(x,\beta) = [\phi(x,\beta,Z_{d(i)}) - \exp(\beta^T Z_{d(i)})] \frac{g_{t_i}(x)(1-x)}{k(t_i)} \tag{35}$$

$$+ \frac{(1-x)\exp(\beta^T Z_{d(i)})}{k(t_i)}(g_{t_i}(x) - k(t_i)) \tag{36}$$

$$+ \exp(\beta^T Z_{d(i)})[(1-x) - (1-x)^{\exp(\beta^T Z_{d(i)})}], \tag{37}$$

where

$$\phi(x,\beta,Z) = \frac{1-(1-x)^{\exp(\beta^T Z)}}{x}.$$

For (35), let $h(x:\beta,Z) = \phi(x,\beta,Z) - \exp(\beta^T Z)$. Then direct calculation yields that

$$m_1 = \sup_{\beta \in B, x \in (0,1/2), \|Z\| \leq M_z} |h'(x:\beta,Z)| < \infty$$

where $h'(x:\beta,Z) = dh(x:\beta,Z)/dx$. Now, since $h(0:\beta,Z) = 0$, the mean value theorem implies that

$$\begin{aligned}
&\sup_{\beta \in B, x \in (0,1/2), \|Z\| \leq M_z} \left|\frac{h(x:\beta,Z)}{x^{\alpha'}}\right| \\
(38) \quad &= \sup_{\beta \in B, x \in (0,1/2), \|Z\| \leq M_z} \left|\frac{h(x:\beta,Z) - h(0:\beta,Z)}{x}\right| x^{1-\alpha'} \\
&= \sup_{\beta \in B, x \in (0,1/2), \|Z\| \leq M_z} |h'(x:\beta,Z)|.
\end{aligned}$$

Hence

$$\sup_{\beta \in B, x \in (0,1/2)} \left\|\frac{(35)}{x^{\alpha'}}\right\| \leq \frac{g^*}{k_*} m_1.$$

Also, it is easy to see that

$$\sup_{\beta \in B, x \in [1/2,1)} \left\|\frac{(35)}{x^{\alpha'}}\right\| \leq D$$

for some constant $D$, since the numerator as well as the denominator is finite.

For (36), conditions (C1) and (C2) imply that

$$m_2 = \sup_{\beta \in B, x \in (0,1), t \in [0,\tau]} \left|\frac{(1-x)(g_t(x) - k(t))}{x^{\alpha'}}\right| < \infty.$$

So

$$\sup_{\beta \in B, x \in (0,1)} \left\|\frac{(36)}{x^{\alpha'}}\right\| \leq \frac{C_B}{k_*} m_2$$

where $C_B = \sup_{\beta \in B, \|Z\| \leq M_z} \exp(\beta^T Z)$.

For (37), let $k(x:\beta,Z) = (1-x) - (1-x)^{\exp(\beta^T Z)}$ and let $k'(x:\beta,Z)$ be the first derivative of $k$ in $x$. Direct calculation yields

$$\sup_{\beta \in B, x \in (0,1/2], \|Z\| \leq M_z} |k'(x:\beta,Z)| < \infty.$$

So we can use a method similar to (38) to show

$$m_3 = \sup_{\beta \in B, x \in (0,1/2], \|Z\| \leq M_z} \left|\frac{k(x:\beta,Z)}{x^{\alpha'}}\right| < \infty.$$

Also, it is true that

$$m_4 = \sup_{\beta \in B, x \in (1/2,1], \|Z\| \leq M_z} \left| \frac{k(x:\beta, Z)}{x^{\alpha'}} \right| < \infty.$$

Hence

$$\left\| \frac{(37)}{x^{\alpha'}} \right\| \leq C_B(m_3 + m_4).$$

Now the proof of (34) for $k = 0$ is done by letting

$$M_0 = m_1 g^*/k_* + D + m_2 C_B/k_* + C_B(m_3 + m_4).$$

The results for $k = 1, 2, 3$ follow from similar arguments. □

PROOF OF LEMMA 1. We can write

$$h_n(\beta) = -\rho_n(\beta) + l_n(\beta) + \sum_{i=1}^{q_n} \log(1 + \zeta_i(\beta)/\chi_i(\beta)),$$

where

$$\chi_i(\beta) = \frac{\exp(\beta^T Z_{d(i)})}{\sum_{j \in R(t_i)} \exp(\beta^T Z_j)}$$

and

$$\zeta_i(\beta) = \int_0^1 \frac{1 - (1-x)^{\exp(\beta^T Z_{d(i)})}}{x} (1-x)^{\sum_{j \in R_n^+(t_i)} \exp(\beta^T Z_j)} \frac{g_{t_i}(x)}{k(t_i)} dx$$
$$- \frac{\exp(\beta^T Z_{d(i)})}{\sum_{j \in R_n(t_i)} \exp(\beta^T Z_j)}.$$

For $\rho_n(\beta)$, using conditions (A2) and (C1), we have

$$(39) \quad \sup_{\beta \in B} \rho_n(\beta) \leq M \sum_{i=1}^n \int_0^{T_i} \int_0^1 (1-x)^{(n-i)c+\varsigma-1} dx \, \lambda(t) \, dt = O(\log n)$$

with probability 1 for some positive constants $M$ and $c$. Similarly, we can show that $\sup_{\beta \in B} \|\rho_n(\beta)^{(k)}\| = O(\log n)$ with probability 1 for $k = 1, 2, 3$.

Let

$$(40) \quad \xi_i(\beta) = \zeta_i(\beta)/\chi_i(\beta).$$

The proof will be complete if we show that $\sup_{\beta \in B, i=1,\ldots,q_n} \|\xi_i^{(k)}(\beta)\| = o(1)$ for $k = 0, 1, 2, 3$. However, we will show $\sup_{\beta \in B, i=1,\ldots,q_n} \|\xi_i^{(k)}(\beta)\| = o(n^{-1/2})$ to use it in the proof of Lemma 2. Since $\sup_{\beta \in B, i=1,\ldots,q_n} \|\chi_i^{(k)}(\beta)\| = O(n^{-1})$

for $k = 0, 1, 2, 3$, it suffices to show $\sup_{\beta \in B, i=1,\ldots,q_n} \|\zeta_i^{(k)}(\beta)\| = o(n^{-3/2})$ for $k = 0, 1, 2, 3$.

For $k = 0$, let $\alpha' = \min\{\alpha, 1\}$ where $\alpha$ is in condition (C2). Since

$$\frac{\exp(\beta^T Z_{d(i)})}{\sum_{j \in R_n(t_i)} \exp(\beta^T Z_j)} = \int_0^1 (1-x)^{\sum_{j \in R_n^+(t_i)} \exp(\beta^T Z_j) - 1}$$

$$\times \exp(\beta^T Z_{d(i)})(1-x)^{\exp(\beta^T Z_{d(i)})} \, dx,$$

we can write

$$\zeta_i(\beta) = \int_0^1 (1-x)^{\sum_{j \in R_n^+(t_i)} \exp(\beta^T Z_j) - 1} \eta_i(x, \beta) \, dx$$

where $\eta$ is defined in (33). Then Lemma A.1 yields

(41)
$$\sup_{\beta \in B, 1 \leq i \leq q_n} \|\zeta_i(\beta)\| \leq \int_0^1 (1-x)^{\sum_{j \in R_n^+(\tau)} \exp(\beta^T Z_j) - 1} x^{\alpha'} M_0 \, dx$$
$$= O(n^{-(\alpha'+1)}) = o(n^{-3/2}),$$

where the last equality is due to the fact that $\alpha' > 1/2$ by condition (C2). Similarly, we can get $\sup_{\beta \in B, i=1,\ldots,q_n} \|\zeta_i^{(k)}(\beta)\| = o(n^{-3/2})$ for $k = 1, 2, 3$. □

PROOF OF LEMMA 2. We have

(42)
$$\|h_n^{(1)}(\hat{\beta})\| \leq \|\rho_n^{(1)}(\hat{\beta})\| + \|l_n^{(1)}(\hat{\beta})\| + \sum_{i=1}^{q_n} \left\| \frac{\xi_i^{(1)}(\hat{\beta})}{1 + \xi_i(\hat{\beta})} \right\|$$

where $\xi_i(\beta)$ is defined in (40). We have shown in the proof of Lemma 1 that $\|\rho_n^{(1)}(\hat{\beta})\| = O(\log n)$ and

$$\sum_{i=1}^{q_n} \left\| \frac{\xi_i^{(1)}(\hat{\beta})}{1 + \xi_i(\hat{\beta})} \right\| = o(\sqrt{n})$$

with probability 1. Since $l_n^{(1)}(\hat{\beta}) = 0$, the proof is done. □

PROOF OF LEMMA 3. Let $\theta_n = (\beta, D_n)$. Let

$$e_i^k(\beta) = \int_0^1 x^k \frac{1 - (1-x)^{\exp(\beta^T Z_{d(i)})}}{x} (1-x)^{\sum_{j \in R^+(t_i)} \exp(\beta^T Z_j)} g_{t_i}(x) \, dx$$

and

$$\tilde{e}_i^k(\beta) = k(t_i) \exp(\beta^T Z_{d(i)}) \frac{\Gamma(\sum_{j \in R(t_i)} \exp(\beta^T Z_j)) \Gamma(k+1)}{\Gamma(\sum_{j \in R(t_i)} \exp(\beta^T Z_j) + k + 1)}.$$

Since
$$\tilde{e}_i^k(\beta) = k(t_i) \int_0^1 x^k \exp(\beta^T Z_{d(i)})(1-x)^{\sum_{j \in R_n(t_i)} \exp(\beta^T Z_j) - 1} \, dx,$$

using Lemma A.1, we have

$$\sup_{\beta \in B} |e_i^k(\beta) - \tilde{e}_i^k(\beta)|$$

$$\leq k^* \sup_{\beta \in B} \left| \int_0^1 x^k (1-x)^{\sum_{j \in R_n^+(t_i)} \exp(\beta^T Z_j) - 1} \eta_i(x, \beta) \, dx \right|$$

$$\leq k^* M_0 \left| \int_0^1 x^{k+\alpha'} (1-x)^{\sum_{j \in R_n^+(t_i)} \exp(\beta^T Z_j) - 1} \, dx \right|$$

$$= O(n^{-(k+\alpha'+1)}) = o(n^{-(k+3/2)})$$

with probability 1.

Now, we can write

$$\sup_{\beta \in B} \left| E(\Delta A^k(t_i) | \theta_n) - \frac{k! \Gamma(\sum_{j \in R(t_i)} \exp(\beta^T Z_j) + 1)}{\Gamma(\sum_{j \in R(t_i)} \exp(\beta^T Z_j) + k + 1)} \right|$$

(43)
$$= \sup_{\beta \in B} \left| \frac{e_i^k(\beta)}{e_i^0(\beta)} - \frac{\tilde{e}_i^k(\beta)}{\tilde{e}_i^0(\beta)} \right|$$

$$\leq \sup_{\beta \in B} \left| \frac{e_i^k(\beta) - \tilde{e}_i^k(\beta)}{e_i^0(\beta)} \right| + \sup_{\beta \in B} \left| \frac{\tilde{e}_i^k(\beta)(e_i^0(\beta) - \tilde{e}_i^0(\beta))}{e_i^0(\beta) \tilde{e}_i^0(\beta)} \right|.$$

Note that $\tilde{e}_i^k(\beta) = O(n^{-(k+1)})$ and hence $e_i^k(\beta) = O(n^{-(k+1)})$. Therefore,

$$(43) = \frac{o(n^{-(k+3/2)})}{O(n^{-1})} + \frac{O(n^{-(k+1)}) o(n^{-3/2})}{O(n^{-2})} = o(n^{-(k+1/2)})$$

with probability 1, and the proof is done. □

DEPARTMENT OF STATISTICS
SEOUL NATIONAL UNIVERSITY
SILLIMDONG KWANAKGU
SEOUL 151-742
REPUBLIC OF KOREA
E-MAIL: ydkim@stats.snu.ac.kr